\documentclass[12pt,reqno]{amsart}
\usepackage{amssymb,delarray}
\usepackage{amsfonts}
\usepackage{epsfig}
\usepackage[all]{xy}

\makeindex{}

\def\le{\leqslant}
\def\ge{\geqslant}

\def\edvo{\rule {6pt}{6pt}}

\begin{document}
\baselineskip=14.pt plus 2pt 

\title[
]{On finite presentability of group $\bf F/[M,N]$}
\author[]{O.V.~Kulikova, A.Yu.~Olshanskii}

\address{\newline O.V.Kulikova \newline Bauman Moscow State Technical University \newline (Moscow, Russia)}
\email{olga.kulikova@mail.ru}

\address{\newline A.Yu.Olshanskii\newline Vanderbilt
University\newline (Nashville, U.S.A.)
\newline Moscow State University
\newline (Moscow, Russia) }
\email{alexander.olshanskiy@vanderbilt.edu}

\dedicatory{} 

\begin{abstract}
The characterization of normal subgroups $M$, $N$ of free group
$F$ for which the quotient group $F/[M,N]$ is finitely presented
is given.

\end{abstract}

\maketitle  \setcounter{tocdepth}{2}
\def\st{{\sf st}}

\setcounter{section}{0}

{\bf 1. Introduction.} Let $F$ be a free group and $M$, $N$ normal
subgroups of $F$. When is the quotient group $F/[M,N]$ finitely
related? In other words: when is the mutual commutator subgroup
$[M,N]$ the normal closure of a finite set of elements?

The case $M=N$ has been already considered. For example, it
follows from Theorem 5.3 of A.L.Shmelkin's paper \cite{shm} that
if $F$ is a finitely generated free group and $N$ is the normal
closure of a finite set of words from $F$, then $F/[N,N]$ is
finitely presented if and only if $N$ has a finite index in $F$.

There is the following lemma in the paper \cite{ar} of
J.Abarbanel, S.Rosset:

{\bf Lemma 1.} {\it If $M$ and $N$ are the normal closures of
finite sets of elements in a finitely generated group $F$ and the
group $MN$ has a finite index in $F$, then $[M,N]$ is the normal
closure of a finite set of elements.}

For completeness, the proof of Lemma 1 will be given in Section 2.

The condition in Lemma 1 that $M$ and $N$ are normally finitely
generated is essential. If $N$ is not finitely generated as normal
subgroup, then $[F,N]$ can be or not be finitely generated as
normal subgroup depending on $N$. In \cite{ar} (referring to
\cite{abels} and \cite{groves}), some examples showing that are
given.

Also in \cite{ar} (and \cite{baum}), the authors consider the
case, when the index of $MN$ in $F$ is infinite. They show that in
this case under some restrictions on $M$ and $N$ the group
$F/[M,N]$ is not finitely related.

The following question remained open: suppose $M$ and $N$ are
non-trivial normal subgroups of the free group $F$. Is it true
that a finite index of $MN$ in $F$ is a necessary condition for
$F/[M,N]$ to be finitely related? The following theorem shows that
the answer to this question is affirmative.

{\bf Theorem 1.} {\it Let $M$ and $N$ be non-trivial normal
subgroups in a free group $F$. Suppose that the group $MN$ has
infinite index in $F$. Then the mutual commutator subgroup $[M,N]$
is not the normal closure of any  finite set of elements.}

So if $M$ and $N$ are the normal closures of finite sets of words
in the finitely generated free group $F$, Lemma 1 and Theorem 1
lead to the following

{\bf Theorem 2.} {\it If $M$ and $N$ are the normal closures of
non-empty finite sets of non-identical words in the finitely
generated free group $F$, then $[M,N]$ is the normal closure of a
finite set of elements if and only if the group $MN$ has a finite
index in $F$.}

{\bf Corollary 1.} {\it Let $F$ be a free group of infinite rank,
and $M$, $N$ non-trivial normal subgroups in $F$. Then the mutual
commutator subgroup $[M,N]$ is not the normal closure of any
finite
set of elements. }\\

{\bf 2. Proofs.}

{\bf Proof of Lemma 1.} The group $L = MN$ is finitely generated
since it has a finite index in the finitely generated group $F$.
Hence there is a finitely generated subgroups $M_1$ and $N_1$ in
$M$ and $N$, respectively, such that $M_1$ and $N_1$ generate $L$,
and $M$ and $N$ are the normal closures of $M_1$ and $N_1$ in $F$.

Let $\{x_i\}$ and  $\{y_j\}$ be finite generating sets of $M_1$
and $N_1$. It is clear that the subgroup $K = [M, N]$ is the
normal closure of all commutators $[x_i^f, y_j]$, where $f\in F$.
So it suffices to show that the relation $[x_i^f, y_j] = 1$ of the
group $F/K$ follows from relations $[x_i^{h_k}, y_j] = 1$ for all
$h_k$ from some set of representatives of left cosets of $L$ in
$F$. Let $f = h_ku_1...u_s$, where $u_1,...,u_s\in \{x_i^{\pm 1}
\}\cup\{y_j^{\pm 1}\}$. The proof will induct on $s$. If $s=0$,
there is nothing to proof. Let $s\geq 1$. If $u_s = y$, where
$y\in \{y_j^{\pm 1}\}$, then the relation $[x_i^f, y_j] = 1$
follows, by assumption, from $x_i^{f'} = y^{-1}x_i^{f'}y =
x_i^{f}$ and $[x_i^{f'}, y_j] = 1$, where $f' = h_ku_1...u_{s-1}$.
If $u_s = x$, where $x\in \{x_i^{\pm 1}\}$, then the relation
$[x_i^f, y_j] = 1$ is equivalent to the relation $[x_i^{f'},
xy_jx^{-1}] = 1$, which follows, by assumption, from $xy_jx^{-1} =
y_j$ and $[x_i^{f'}, y_j] = 1$. \edvo\vspace{3mm}

{\bf Proof of Theorem 1.} Below we use the following notation:

$X$, a set of free generators of the free group $F$;

$L$, the product $MN$ of the non-trivial groups $M$ and $N$;

$\gamma_i(L)$, the $i$-th term of the lower center series for $L$;

$G$, the quotient group $F/L$.

It is well-known that the intersection of all terms of the lower
center series $\gamma_i(L)$ in the free group is trivial.
Therefore there is a number $c\ge 1$ such that the group $M$ is
contained in $\gamma_c(L)$, but it is not contained in
$\gamma_{c+1}(L)$. Without loss of generality we can assume that
$N$ is not contained in $\gamma_2(L)$, since $L = MN$. Then
$[M,N]$ is contained in $\gamma_{c+1}(L)$.

To prove Theorem 1, it suffices to show that $[M,N]$ is not
contained in the normal closure of any finite subset $Y$ from
$\gamma_{c+1}(L)$. To show this, we use verbal wreath products,
introduced by A.L. Shmelkin in \cite{shm}. More precisely, we use
the $c$-nilpotent wreath product $S$ of a free $c$-nilpotent group
$T(X')$ and $G = F/L$, where $X'$ is of the same cardinality as
$X$.

By definition, $S$ is a semidirect product of a group $W$ and $G$,
where $W$ is the $c$-nilpotent product of  $|G|$ copies of the
group $T(X')$. Moreover, $W$ is the free $c$-nilpotent group with
basis $\{a_i| i=(x,h), x\in X', h\in G\}$. The group $G$ acts on
$W$ by the rule $g^{-1}a_ig = a_j$, where $i=(x,h)$, $j = (x,hg)$.

By Theorem 2.1 \cite{shm} (see also Addition 1 to \cite{neim}),
the homomorphism $\psi$ from $F$ to $S$, mapping $x\in X$ to
$\tilde x a_{(x,1)}$, where $\tilde x = xL\in G = F/L$, has the
kernel $\gamma_{c+1}(L)$, that is, the image $H$ of $F$ is
isomorphic to the quotient group $F/\gamma_{c+1}(L)$. According to
the properties of this homomorphism, the intersection of $H$ and
$W$ is the image of the subgroup $L$, and $HW = S$.\vspace{3mm}

Further we consider an arbitrary subset $Z$ in $G$ such that $Z =
Z^{-1}$.

We say that two bi-indexes $i = (x',h)$ and $j = (x'',g)$, where
$x', x''\in X'$, $h, g\in G$, are {\it $Z$-close}, if the quotient
$hg^{-1}$ of their second components belongs to $Z$, and are {\it
$Z$-distant} otherwise.

For a given subset $Z$, let us construct  an auxiliary semidirect
product $S_Z$ of groups $W_Z$ and $G$ (so that $S_Z = S$ for $Z =
G$). The group $W_Z$ is defined as follows. The generators for
$W_Z$ are the same $a_i$ as in $W$ with the same bi-indexes $i$.
The defining relations for $W_Z$ are left-normalized commutators
$[a_{i_1},..., a_{i_{c+1}}]$ of generators such that any two
bi-indexes $i_k$ and $i_l$ in such commutator are $Z$-close. It is
clear that these defining relations are invariant under the action
of the group $G$ on $W_Z$ by conjugation, which shifts the second
components of bi-indexes of the generators $a_i$ by the rule
$g^{-1}a_ig = a_j$, where $i = (x,h)$, $j = (x,hg)$.

By such construction, if $Z = \emptyset$, then $W_\emptyset$ is an
absolutely free group, and $S_\emptyset$ is the semidirect product
of this free group $W_\emptyset$ and $G$. If $Z = G$, then $W_G =
W$ and $S_G = S$. The other cases are intermediate between
$W_\emptyset$ and $W$ and between $S_\emptyset$ and $S$,
respectively.

Moreover, if $Z\subset Z'\subseteq G$, the group $S_Z$ is
naturally mapped onto $S_{Z'}$ by the homomorphism $\phi_{Z\subset
Z'}$ identical on $G$ and the set of generators $\{a_i\}$.

For any subset $Z\subseteq G$, one can define the homomorphism
$\psi_Z$ from the free group $F$ to $S_Z$, mapping $x\in X$ to
$\tilde x a_{(x,1)}$, where $\tilde x = xL\in G = F/L$ and
$a_{(x,1)}$ are generators of $W_Z$. Since $S_Z = GW_Z$ and
$\tilde x$ generate $G$, we have that $S_Z = H_ZW_Z$, where $H_Z =
\psi_Z(F)$. Besides, $\psi = \psi_G$ and $\psi = \phi_{Z\subseteq
G}\psi_Z$ for any $Z\subseteq G$.\vspace{3mm}

Recall that $\gamma_{c+1}(L)$ is the kernel of the homomorphism
$\psi: F\rightarrow S$. On the other hand, this kernel is the
union of the kernels of the homomorphisms $\psi_Z$ to $S_Z$ for
all finite subsets $Z\subset G$. Hence any finite subset $Y$ from
$\gamma_{c+1}(L)$ belongs to the kernel of a homomorphism $\psi_Z:
F\rightarrow S_Z$ for some finite subset $Z$.

Therefore, to show that $[M,N]$ is not contained in the normal
closure of any finite subset $Y$ from $\gamma_{c+1}(L)$, it
suffices to  show that, for any finite subset $Z\subseteq G$,
$[M,N]$ is not mapped to the trivial subgroup by the homomorphism
$\psi_Z$ of the free group $F$ to $S_Z$. Below we use the
following simple lemma.

{\bf Lemma 2.} {\it Let $u$ and $v$ be not identical elements of
$W_Z$ such that the bi-index of any generator in some writing of
the element $u$ is $Z$-distant to the bi-index of any generator in
some writing of the element $v$. Then $w^{-1}uw$ and $v$ do not
commute for any element $w\in W_Z$. }

{\bf Proof of Lemma 2.} Let $A$ be a subgroup of $W_Z$, generated
by such $a_i$ that occur in the writing of $u$. Similarly, the
subgroup $B$ is defined for $v$. Consider the free product $A*B$.
Let us define a homomorphism from the group $W_Z$ to $A*B$ as
follows. It leaves fixed generators of $A$ and $B$ and maps the
other generators to the identity. This definition is well-defined,
since the defined homomorphism remains the defining relations of
$W_Z$. In addition the image of $u$ belongs to $A$ and the image
of $v$ belongs to $B$. Now Lemma 2 follows from the well-known
fact that if not identical elements are conjugated with elements
from different factors of the free product, then they do not
commute. \edvo

Now we complete the proof of Theorem 1.

Since the groups $M$ and $N$ do not belong to $\gamma_{c+1}(L)$,
there exist elements $u\in M$ and $v\in N$, whose images in the
group $S$ are not equal to the identity. Moreover, their images
under the homomorphism $\psi$ belong to $W$, since $M$ and $N$
belong to $L$. Therefore, their images $u'$ and $v'$ under the
homomorphism $\psi_Z$ from $F$ to $S_Z$ are not equal to the
identity in the group $W_Z$, since $\psi = \phi_{Z\subseteq
G}\psi_Z$.

Recall that the second components of bi-indexes of generators
$a_i$ of the group $W_Z$ are multiplied from the right by $g$
under the conjugation by an element $g\in G$. Since $Z$ is finite
and $G$ is infinite, there exists an element $g\in G$ such that
the bi-index of any generator in some writing of the following
element
\begin{equation} u'' = g^{-1}u'g\label{eq1}
\end{equation} is $Z$-distant from the bi-index of any generator in some writing of the element $v'$.

Since $G\le S_Z = H_ZW_Z$, for the element $g$, there exist $h\in
H_Z$ and $w\in W_Z$ such that $g = hw^{-1}$. According to
(\ref{eq1}), $u'' = wh^{-1}u'hw^{-1}$. Hence, $w^{-1}u''w =
h^{-1}u'h$. Since $h\in H_Z=\psi_Z(F)$, $u'\in\psi_Z(M)$, we have
that $h^{-1}u'h\in\psi_Z(M)$. Therefore, $w^{-1}u''w\in\psi_Z(M)$.

Since the bi-index of any generator in the writing of the element
$u''$ is $Z$-distant from the bi-index of any generator in the
writing of the element $v'$, it follows from Lemma 2 that the
elements $w^{-1}u''w$ and $v'$ do not commute. Hence, taking into
account that $w^{-1}u''w\in\psi_Z(M)$ and $v'\in\psi_Z(N)$, we
have that the image of $[M,N]$ is non-trivial in $S_Z$.
\edvo\vspace{3mm}

{\bf Proof of Corollary 1.} If $F$ has infinite rank and $M$, $N$
are finitely generated as normal subgroups, then $MN$ has infinite
index in $F$. Consequently, Corollary 1 follows from Theorem 1.

If $M$ and (or) $N$ are infinitely generated as normal subgroups,
then the group $M$ (respectively, $N$) can be presented as the
union of subgroups $M_k$ (respectively, $N_k$), where $M_k$
($N_k$) is finitely generated as normal subgroup. Similar to the
previous paragraph, we have that the group $[M_k,N_k]$ is not
generated by any finite set of elements as normal subgroup. Hence,
the union of all $[M_k,N_k]$ can not be the normal closure of any
finite set of elements. Since the union of all subgroups
$[M_k,N_k]$ is equal to the group $[M,N]$, Corollary 1 is proved.
\edvo\vspace{3mm}

{\bf Remark.} Let $F$ be a free group, and $M$, $N$ non-trivial
normal subgroups in $F$ such that the group $MN$ has infinite
index in $F$ or $F$ has infinite rank. Suppose that $[M,N] = M\cap
N$. For example, it is true when the union $S\cup T$ satisfies the
small cancellation condition $C'(\lambda)$ ($\lambda \leq 1/6)$,
where $S$ (respectively, $T$) is a non-empty symmetrized set of
cyclically reduced words from $F$ such that $S$ (respectively,
$T$) generates $M$ (respectively, $N$) as normal subgroup (see,
for example, \cite{ratcl, ok}). Then it follows from Theorem 1
(Corollary 1) that $M\cap N$ is not the normal closure of any
finite set of elements, even if the groups $M$ and $N$ are
finitely generated as normal subgroups. That is contrasting to the
well-known A.G.Howson's Theorem (see, for example, Proposition
3.13 in \cite{lind}): if subgroups $H$ and $K$ in a free group are
finitely generated, then their intersection $H\cap K$ is finitely
generated.

The paper is partially supported by RFBR 05-01-00895, the research
of the second author is also supported by NSF DMS-0245600 and
DMS-04556881.



\begin{thebibliography}{9}


\bibitem {ar} {\it Abarbanel J.,
Rosset S.} The Schur multiplier of $F/[R,S]$ // J. Pure and Appl.
Algebra. 2005. {\bf 198}.  1-8.

\bibitem {abels} {\it Abels H.}
An Example of a Finitely Presented Solvable Group // Homological
Group Theory. London Math. Soc. Lec. Notes Ser. Cambridge:
Cambridge University Press. 1979. Vol. {\bf 36}. 205-211.

\bibitem {baum} {\it
Baumslag G., Strebel R., Thomson W.}  On the multiplicator of
$F/\gamma_cR$ // J. Pure and Appl. Algebra. 1980. {\bf 16}.
121-132.

\bibitem {groves} {\it Groves J.R.J.} Finitely presented
centre-by-metabelian groups // J. London Math.Soc. 1978. {\bf 18},
N 2. 65-69.

\bibitem {ratcl} {{\it
Guti\'{e}rrez M.A., Ratcliffe J.G.}  On the second homotopy group
// Quart. J. Math. Oxford (2). 1981. {\bf 32}. 45-55.}

\bibitem {ok} {{\it Kulikova O.V.}  On intersections of normal
subgroups in free groups // Algebra and Discrete Mathematics.
2003. N 1. 36-67.}

\bibitem {lind} { {\it Lindon R.S., Schupp P.E.}
Combinatorial group theory // Berlin; Heidelberg; N.Y.:
Springer-Verlag. 1977. }

\bibitem {neim1} {\it Neumann H.}  Varieties of Groups // Springer-Verlag: Berlin-Heidelberg-New York, 1967.

\bibitem {neim} {\it Neumann H.}   Varieties of Groups // Moscow: Mir,
1969 (in russian).

\bibitem {shm} {\it Shmelkin A.L.} Wreath products and varieties of groups // Izv. Akad. Nauk SSSR. Ser. Mat.,
{\bf 29} N 1, 1965, 149-170 (in russian).



\end{thebibliography}
\end{document}